\documentclass[11pt]{article}
\usepackage[utf8]{inputenc}
\usepackage{amsmath, amssymb, amsthm}
\usepackage[margin=1in]{geometry}
\usepackage{color}
\usepackage{hyperref}
\usepackage[english]{babel}
\usepackage[style=numeric, backend=biber]{biblatex}
\newtheorem{theorem}{Theorem}

\newtheorem{definition}[theorem]{Definition}
\newtheorem{remark}[theorem]{Remark}
\newcommand{\R}{\mathbb{R}}
\title{Strict Lyapunov function for the one-dimensional linear wave equation with a locally distributed damping}
\author{Jean-Michel Coron\thanks{Laboratoire Jacques-Louis Lions, Sorbonne Universit\'{e}, Universit\'{e} Paris Cit\'{e}, CNRS, INRIA, \'{e}quipe CAGE, Paris, France. (jean-michel.coron@sorbonne-universite.fr)}}
\date{}

\begin{document}
\maketitle
\begin{abstract}
We construct an explicit strict Lyapunov function for the one-dimensional linear wave equation with a locally distributed damping and Dirichlet boundary conditions.
\end{abstract}
\noindent \textbf{Keywords:} Linear one-dimensional wave equation; Strict Lyapunov function; locally distributed damping.

\noindent \textbf{2020 Mathematics Subject Classification:} 93D30; 35L05; 93C20.
\section{Introduction and assumptions} We consider the damped wave equation on $[0,L]$:
\begin{equation} \label{eq-wave}
\left\{
\begin{array}{l} y_{tt} - y_{xx} + a(x)y_t = 0,\quad t \geq 0,\quad x \in (0,L),
\\ y(t,0) = y(t,L) = 0,\quad t \geq 0.
\end{array} \right.
\end{equation}
We assume that
\begin{gather}
\label{a-continuous-and-geq0}
a \in L^\infty(0,L),
\\
\label{anot0} a \geq 0 \quad \text{almost everywhere},
\end{gather}
and that there exist $0 \leq x_0 < x_1 \leq L$ and $a_{\text{min}} \in (0,+\infty)$
such that
\begin{equation}
\label{aona0a1} a(x) \geq a_{\text{min}} > 0,
\quad \text{for almost every } x \in (x_0,x_1).
\end{equation}
The energy of the system is defined by
\begin{equation}
\label{def:E}
E(y, y_t) = \frac{1}{2} \int_0^L \left( y_t^2 + y_x^2 \right)\, dx.
\end{equation}
Taking the time derivative of $E$ along solutions of \eqref{eq-wave}
yields the classical dissipation law (see also \eqref{dotE-new} below)
\begin{equation}
\label{eq:dot_E}
\dot E = - \int_0^L a(x) y_t^2\, dx \leq 0.
\end{equation}
Since $a(x)$ may vanish, even on a nonempty open interval, energy dissipation alone does
not guarantee global control of $y_t^2$,
nor does it provide any control of $y_x^2$.
Nevertheless, it is already known (see Theorem~\ref{th:exp-stability} below)
that $(0,0)\in H^1_0(0,L)\times L^2(0,L)$
is exponentially stable with respect to the $H^1_0(0,L)\times L^2(0,L)$ norm.
Our goal is to construct a strict Lyapunov function by adding suitable terms to $E$.
This will allow us to recover exponential stability directly
and may also open the way to studying robustness
issues in a simple and quantitative manner.

The paper is organized as follows.
\begin{enumerate}
\item In Section~\ref{sec:wellposed},
we recall classical results on the well-posedness of the Cauchy problem associated with \eqref{eq-wave},
as well as classical exponential stability results for this equation.
\item In Section~\ref{sec:constructionLyap},
we introduce the terms that are added to $E$ in order to obtain a
strict Lyapunov function in the framework of weak solutions. Some of these terms are local
(Section~\ref{sec-local-terms}) and of a classical form, whereas our construction
also requires a nonlocal term (Section~\ref{sec:nonlocal}) which is less classical.
\item In Section~\ref{sec:strict-for-strong}, we construct a strict Lyapunov function
in the framework of strong solutions. In contrast to the previous Lyapunov function, this one contains only local terms.
\item Finally, Section~\ref{sec:conclusion-and-perspectives} provides the conclusion together with some perspectives.
\end{enumerate}
\begin{enumerate}
\item Transforming a nonstrict Lyapunov function into a strict one is a well-explored problem in finite-dimensional control theory, often referred to as strictification \cite{2009-Malisoff-Mazenc-book}. Despite the availability of powerful methods in that setting, their adaptation to infinite-dimensional systems, such as the wave equation \eqref{eq-wave} and broader classes of partial differential equations, presents a significant theoretical challenge.
\item When distributed damping is replaced by classical boundary damping, the wave equation can be transformed into a $2\times 2$ linear one-dimensional hyperbolic system. This allows for the construction of an explicit strict Lyapunov function; see, for example, \cite{2008-Coron-Bastin-d-Andrea-Novel-SICON,2016-Bastin-Coron-book,2016-Prieur-Tarbouriech-Gomes-da-Silva-IEEE-TAC}. However, this is no longer possible in the presence of distributed damping. Such damping introduces a source term into the $2\times 2$ linear one-dimensional hyperbolic system, which prevents the construction of a strict Lyapunov function using known methods.
\end{enumerate}

\section{Well-posedness of the Cauchy problem and exponential stability}
\label{sec:wellposed}
In this section, we recall classical results on the well-posedness of the Cauchy problem associated with our wave equation
\eqref{eq-wave}, namely, the solution to
\begin{equation}
\label{Cauchy-wave}
\left\{
\begin{array}{l}
y_{tt} - y_{xx} + a(x)y_t = 0, \quad t \geq 0, \quad x \in (0,L), \\
y(t,0) = y(t,L) = 0, \quad t \geq 0, \\
y(0)=y_0 \text{ and } y_t(0)=y_1,
\end{array}
\right.
\end{equation}
where $y_0:(0,L)\rightarrow \R$ and $y_1:(0,L)\rightarrow \R$ are given. We also recall exponential stability results. In \eqref{Cauchy-wave} and in what follows, we use the standard notation $y(t)=y(t,\cdot)$.

We use the energy method framework. Another classical alternative is the semigroup framework and the use of the Hille--Yosida and Lumer--Phillips theorems; see, for example, \cite{1983-Pazy-book}, \cite{1994-Komornik-book}, and \cite{2024-Chitour-Nguyen-Roman-A}.

We start with a classical definition of a weak solution to the Cauchy problem \eqref{Cauchy-wave}.
\begin{definition}[Weak solution]
Let $y_0 \in H^1_0(0,L)$ and $y_1 \in L^2(0,L)$. A function $y$ is called a weak solution of \eqref{Cauchy-wave} if:
\begin{gather}
y \in C^0([0,+\infty); H^1_0(0,L))\cap  C^1([0,+\infty);L^2(0,L)),
\\
y_{tt} \in C^0([0,+\infty); H^{-1}(0,L)),
\\
y(0)=y_0 \text{ and } y_t(0)=y_1,
\end{gather}
and for every test function $v \in H^1_0(0,L)$ and for every $t \in [0,+\infty)$, the following identity holds:
\begin{equation}
\langle y_{tt}(t), v \rangle_{H^{-1},H^1_0} + \int_0^L y_x(t) v_x \, dx + \int_0^L a(x) y_t(t) v \, dx = 0.
\end{equation}
\end{definition}

Based on this definition, we have the following theorem.
\begin{theorem}[Existence, uniqueness, continuity, and strong solutions]
\label{th:Cauchy-pb}
Given $y_0 \in H^1_0(0,L)$ and $y_1 \in L^2(0,L)$, there exists a unique weak solution $y$ to problem \eqref{Cauchy-wave}, and for every $T>0$, the map $(y_0,y_1)\in H^1_0(0,L) \times L^2(0,L) \rightarrow y \in C^0([0,T]; H^1_0(0,L))\cap  C^1([0,T]; L^2(0,L))$ is continuous.
Moreover, if $y_0 \in H^2(0,L)\cap H^1_0(0,L)$ and $y_1 \in H^1_0(0,L)$,
\begin{equation}
\label{eq-sol-reg}
y \in C^0([0,+\infty); H^2(0,L) \cap H^1_0(0,L)) \cap C^1([0,+\infty); H^1_0(0,L)) \cap C^2([0,+\infty); L^2(0,L)).
\end{equation}
\end{theorem}
Solutions $y$ satisfying \eqref{eq-sol-reg} are called strong solutions. The proof of Theorem~\ref{th:Cauchy-pb}, which is now classical, is omitted. It relies on the Galerkin approximation method and energy estimates. See, for example, \cite[Chapitre 3, Section 8]{1968-Lions-Magenes-book} and \cite[Section 7.2]{1998-Evans-book}. 
For the less-studied case of the $L^p$ setting and a different approach relying on d'Alembert's formula, see \cite{2024-Chitour-Nguyen-COCV}.

The exponential stability statements in the following theorem are found in \cite[Theorem 3]{1975-Rauch-Taylor-IUMJ}.
\begin{theorem}
\label{th:exp-stability}
The wave equation \eqref{eq-wave} is exponentially stable in both spaces $H^1_0(0,L)\times L^2(0,L)$ and $\left(H^2(0,L)\cap H^1_0(0,L)\right)\times H^1_0(0,L)$: there exist $C>0$ and $\gamma>0$ such that,
for the solution $y$ of the Cauchy problem \eqref{Cauchy-wave},
\begin{gather}
\label{exp-weak}
\|y(t)\|_{H^1}+\|y_t(t)\|_{L^2}\leq Ce^{-\gamma t } \left(\|y_0\|_{H^1}+\|y_1\|_{L^2}\right),\quad
\forall t\geq 0,\, \forall y_0\in H^1_0(0,L), \, \forall y_1 \in L^2(0,L),
\end{gather}
and
\begin{multline}
\label{exp-strong}
\|y(t)\|_{H^2}+\|y_t(t)\|_{H^1}\leq \\
Ce^{-\gamma t } \left(\|y_0\|_{H^2}+\|y_1\|_{H^1}\right), \quad
\forall t\geq 0,\, \forall y_0\in H^2(0,L)\cap H^1_0(0,L), \, \forall y_1 \in H^1_0(0,L).
\end{multline}
\end{theorem}
A key ingredient in the proof of \eqref{exp-weak} is the following observability inequality. Let $T\geq 2L$.
Then, there exists $C>0$ such that
\begin{equation}
\label{observability-ineq}
E(y(T),y_t(T)) \leq C \int_0^T \int_0^L a(x)y_t^2 \, dx \, dt;
\end{equation}
see \cite[(9), p.~86]{1975-Rauch-Taylor-IUMJ}.

\begin{remark}
For more complex  boundary conditions, see, for example, \cite{2016-Prieur-Tarbouriech-Gomes-da-Silva-IEEE-TAC, 2024-Chitour-Nguyen-Roman-A, 2026-Dahmani-Chitour-Nguyen-Roman-JDE}. The pioneering articles on the multi-dimensional case are \cite{1986-Ho-CRAS}, \cite[Chapitre 7]{1988-Lions-book-vol-1}, and \cite{1992-Bardos-Lebeau-Rauch-SICON}. References \cite{1986-Ho-CRAS} and \cite[Chapitre 7]{1988-Lions-book-vol-1} rely on the multiplier method (see Remark~\ref{rem:cross-terms} below), while \cite{1992-Bardos-Lebeau-Rauch-SICON} relies on microlocal analysis and provides an (almost) optimal condition on the support of $a$ for exponential stability. See also \cite{2016-Laurent-Leautaud-COCV} for quantitative estimates. The case of nonlinear wave equations has been studied in numerous articles. See, for example, \cite{1988-Haraux-Zuazua-ARMA, 1990-Komornik-Zuazua-JMPA, 1991-Zuazua-JMPA, 2023-Vanspranghe-Ferrante-Prieur, 2026-Prieur-Tarbouriech-A}, which rely on the multiplier method, and \cite{2004-Alabau-CRAS, 2005-Alabau-AMO, 2010-Alabau-JDE}, which rely on new convexity properties of an explicit function connected to energy dissipation.
\end{remark}

\section{Construction of a strict Lyapunov function for weak solutions}
\label{sec:constructionLyap}
We start by recalling the proof of \eqref{eq:dot_E}, allowing us to introduce some notation and density arguments used later on for other functionals. Let $E: (f,g)\in H^1_0(0,L)\times L^2(0,L) \rightarrow E(f,g)\in \R$ be defined by
\begin{gather}
\label{def:E-f-g}
E(f,g):=\frac{1}{2}\int_0^L f_x^2+g^2\,dx.
\end{gather}
Let $(y_0,y_1) \in H^1_0(0,L) \times L^2(0,L)$. Let $y$ be the unique weak solution to the Cauchy problem \eqref{Cauchy-wave}. With a slight abuse of notation, let $E:t\in [0,+\infty)\rightarrow E(t)\in \R$ be defined by
\begin{equation}
\label{def:E(t)}
E(t):=E(y(t),y_t(t))=\frac{1}{2}\int_0^Ly_x(t)^2+y_t(t)^2 \, dx.
\end{equation}
If
\begin{equation}
(y_0,y_1) \in \left(H^2(0,L)\cap H^1_0(0,L)\right)\times H^1_0(0,L),
\label{reg-y0-y1}
\end{equation}
then, by Theorem~\ref{th:Cauchy-pb}, one has \eqref{eq-sol-reg} and therefore $E\in C^1([0,+\infty))$,  and
\begin{equation}
\label{eq:dot_E-int}
\dot E = \int_0^L (y_x y_{xt} + y_t y_{tt})\, dx
= \int_0^L y_t (y_{xx} - a(x)y_t)\, dx + \int_0^L y_x y_{xt}\, dx.
\end{equation}
Hence, using an integration by parts and the fact that $y(t)$ vanishes at $x=0$ and $x=L$, we get that
\begin{equation}
\label{dotE-new}
\dot E = - \int_0^L a(x) y_t^2\, dx.
\end{equation}
Note that, using the density of $\left(H^2(0,L)\cap H^1_0(0,L)\right)\times H^1_0(0,L)$ in $H^1_0(0,L)\times L^2(0,L)$ and the continuous dependence of $y$ on the initial condition $(y_0,y_1)$ stated in Theorem~\ref{th:Cauchy-pb}, we deduce that if $(y_0,y_1)$ is only in $H^1_0(0,L)\times L^2(0,L)$, then $E$ remains of class $C^1$ and \eqref{dotE-new} still holds.

We continue by adding local terms to $E$ (Section \ref{sec-local-terms}). It turns out that we were not able to obtain a strict Lyapunov function by adding only local terms. For this reason, we add a nonlocal term in Section~\ref{sec:nonlocal}.

\subsection{Local terms}
\label{sec-local-terms}
Let $p$ and $q$ be in $C^2([0,L])$, and let us define $ M: (f,g)\in H^1_0(0,L)\times L^2(0,L) \rightarrow M(f,g)\in \R$ by
\begin{equation}
\label{def:M}
M(f, g) := \int_0^L p(x) f g \, dx + \int_0^L q(x) g f_x\, dx .
\end{equation}

We select $q$ so that
\begin{gather}
\label{q0>0}
q(0) \geq 0,
\\
\label{qL<0}
q(L) \leq 0,
\\
\label{prop:deltaq}
\text{there exists $\delta_q>0$ such that } q'(x) \geq \delta_q > 0, \quad \forall x \in [0,x_0]\cup [x_1,L].
\end{gather}

Once $q$ is fixed, we select $p$ so that
\begin{gather}
p(x) = 0, \quad \forall x \in [0,x_0]\cup [x_1,L],
\label{support-p}
\\
\text{there exists $\delta_p\in (0,\delta_q /2)$ such that } p(x) + \frac{1}{2}q'(x) \geq \delta_p > 0, \quad \forall x \in [x_0,x_1].
\label{prop-deltap}
\end{gather}
Let us check the existence of such $p$. Let $p_0\in C^2([0,L];[0,+\infty))$ be such that
\begin{gather}
p_0=0 \text{ in }[0,x_0]\cup [x_1,L],
\\
p_0>0 \text{ in } (x_0,x_1),
\end{gather}
and then, using \eqref{prop:deltaq}, it suffices to define $p=\lambda p_0$ with $\lambda \geq 1$ large enough.

With a slight abuse of notation already used for $E$, let
\begin{equation}
M(t):=M(y(t),y_t(t)), \quad t\geq 0.
\end{equation}
Let us compute and estimate the time derivative of $M$. As for $E$, we  start with the case where \eqref{reg-y0-y1} holds. Then, according again to Theorem~\ref{th:Cauchy-pb}, $y$ has the regularity given by~\eqref{eq-sol-reg}. Hence, $M$ is of class $C^1$, and we have
\begin{align}
\dot M &= \int_0^L q (y_{tt} y_x + y_t y_{xt}) \, dx + \int_0^L p(y_t^2 + y y_{tt}) \, dx
\nonumber \\
&= \int_0^L q \left( (y_{xx} - a y_t)y_x + \frac{1}{2}(y_t^2)_x \right) \, dx + \int_0^L p y_t^2 \, dx + \int_0^L p y(y_{xx} - a y_t) \, dx.
\end{align}
Using integration by parts:
\begin{align}
\dot M &= \frac{1}{2} \left[ q(x) \left( y_x^2 + y_t^2 \right) \right]_0^L - \frac{1}{2} \int_0^L q'(x) (y_x^2 + y_t^2) \, dx - \int_0^L a(x) q(x) y_t y_x \, dx \nonumber \\
&\quad + \int_0^L p(x) y_t^2 \, dx - \int_0^L p(x) y_x^2 \, dx - \int_0^L p'(x) y y_x \, dx - \int_0^L a(x) p(x) y y_t \, dx.
\label{eq:dot_Mnb1}
\end{align}
Since $y_t=0$ on the boundary, the boundary term is
\begin{equation}
\label{boundary-term-dotM}
\frac{1}{2} q(L)y_x(L)^2 - \frac{1}{2} q(0)y_x(0)^2,
\end{equation}
which is $\leq 0$ by \eqref{q0>0} and \eqref{qL<0}. Hence, rearranging the quadratic terms $y_x^2$ and $y_t^2$ in \eqref{eq:dot_Mnb1}, one gets:
\begin{multline}
\dot M \leq - \int_0^L \left( p(x) + \frac{1}{2}q'(x) \right) y_x^2 \, dx + \int_0^L \left( p(x) - \frac{1}{2}q'(x) \right) y_t^2 \, dx \\
- \int_0^L p'(x) y y_x \, dx - \int_0^L a(x) q(x) y_t y_x \, dx - \int_0^L a(x) p(x) y y_t \, dx.
\label{dotMleqnb1}
\end{multline}

Note that, using the density of $\left(H^2(0,L)\cap H^1_0(0,L)\right)\times H^1_0(0,L)$ in $H^1_0(0,L)\times L^2(0,L)$ and the continuous dependence of $y$ on the initial condition $(y_0,y_1)$ stated in Theorem~\ref{th:Cauchy-pb}, we deduce that \eqref{dotMleqnb1} still holds in the distributional sense (with $M$ being a priori only continuous) if $(y_0,y_1)$ is only in $H^1_0(0,L)\times L^2(0,L)$, which we now assume:
\begin{multline}
\dot M \leq - \int_0^L \left( p(x) + \frac{1}{2}q'(x) \right) y_x^2 \, dx + \int_0^L \left( p(x) - \frac{1}{2}q'(x) \right) y_t^2 \, dx \\
- \int_0^L p'(x) y y_x \, dx - \int_0^L a(x) q(x) y_t y_x \, dx - \int_0^L a(x) p(x) y y_t \, dx
 \text{ in } \mathcal{D}'(0,+\infty),
\label{dotMleqnb1-dist}
\end{multline}
where $\mathcal{D}'(0,+\infty)$ is the set of distributions on $(0,+\infty)$.

\begin{remark} If, instead of \eqref{q0>0} and \eqref{qL<0}, one requires
\begin{equation}
\label{qbord=0}
q(0)=q(L)=0,
\end{equation}
then $M$ is of class $C^1$ on $[0,+\infty)$, and \eqref{dotMleqnb1-dist}, as well as \eqref{dotMleqnb2}, \eqref{dotMleqnb3}, \eqref{dotV1leq}, \eqref{dotV1leqnb2}, \eqref{eq:V1_bound}, \eqref{dotVleq-gammaV}, and \eqref{dotV1leq-new} below, are inequalities between two continuous functions in the classical sense. Moreover, throughout this article, \eqref{q0>0} and \eqref{qL<0} can be replaced by \eqref{qbord=0}. However, it is advantageous to allow
\begin{equation}
\label{qbordsigned}
q(0) >0 \text{ and } q(L)<0
\end{equation}
since it provides robustness with respect to the actual behavior of $y$ on the boundary of $[0,L]$. For this reason, we prefer to keep \eqref{q0>0} and \eqref{qL<0} instead of \eqref{qbord=0}.
\end{remark}

We bound the cross terms using the Young inequality. Let $\eta > 0$. Using also \eqref{support-p} for \eqref{cross-2}, we have:
\begin{gather}
\left| \int_0^L a q y_t y_x \, dx \right| \leq \int_0^L \eta y_x^2 \, dx + \int_0^L \frac{a(x)^2 q(x)^2}{4\eta} y_t^2 \, dx \leq \eta \int_0^L y_x^2 \, dx + \frac{\|a\|_{L^\infty} \|q\|_{L^\infty}^2}{4\eta} \int_0^L a(x) y_t^2 \, dx,
\label{cross-1}
\end{gather}
\begin{multline}
\left| \int_0^L a p y y_t \, dx \right| \leq \int_{x_0}^{x_1} \frac{1}{2} y^2 \, dx + \int_0^L \frac{a(x)^2 p(x)^2}{2} y_t^2 \, dx \leq \frac{1}{2} \int_{x_0}^{x_1}y^2 \, dx
\\
+ \frac{\|a\|_{L^\infty} \|p\|_{L^\infty}^2}{2} \int_0^L a(x) y_t^2 \, dx.
\label{cross-2}
\end{multline}
Moreover, using integration by parts and \eqref{support-p}, we have
\begin{equation}
\label{cross-3}
-\int_0^L p'(x) y y_x \, dx = - \int_0^L p'(x) \frac{1}{2} (y^2)_x \, dx = \frac{1}{2} \int_0^L p''(x) y^2 \, dx \leq \frac{1}{2} \|p''\|_{L^\infty} \int_{x_0}^{x_1} y^2 \, dx.
\end{equation}
Substituting these three bounds \eqref{cross-1}, \eqref{cross-2}, and \eqref{cross-3} into the inequality \eqref{dotMleqnb1-dist} for $\dot M$ and using \eqref{prop-deltap}, we obtain:
\begin{multline}
\dot M \leq - (\delta_p - \eta) \int_0^L y_x^2 \, dx
\\
+ \int_0^L \left( p(x) - \frac{1}{2}q'(x) \right) y_t^2 \, dx
+ C_1 \int_{x_0}^{x_1} y^2 \, dx + C_2(\eta) \int_0^L a(x) y_t^2 \, dx \text{ in } \mathcal{D}'(0,+\infty),
\label{dotMleqnb2}
\end{multline}
where
\begin{equation}
C_1 := \frac{1}{2}(1 + \|p''\|_{L^\infty}) \text{ and } C_2(\eta) := \|a\|_{L^\infty} \left( \frac{\|q\|_{L^\infty}^2}{4\eta} + \frac{\|p\|_{L^\infty}^2}{2} \right).
\end{equation}
We choose
\begin{equation}
\label{choice-eta}
\eta := \frac{\delta_p}{2}.
\end{equation}
Then, \eqref{dotMleqnb2} becomes
\begin{multline}
\dot M \leq - \frac{\delta_p}{2} \int_0^L y_x^2 \, dx + \int_0^L \left( p(x) - \frac{1}{2}q'(x) \right) y_t^2 \, dx
\\
+ C_1 \int_{x_0}^{x_1} y^2 \, dx + C_2(\delta_p/2) \int_0^L a(x) y_t^2 \, dx \text{ in } \mathcal{D}'(0,+\infty).
\label{dotMleqnb3}
\end{multline}

We combine the energy $E$ defined in \eqref{def:E-f-g} and the local multiplier $M$ defined in \eqref{def:M}: we define
\begin{equation}
\label{defV1}
V_1:= E + \varepsilon_1 M,
\end{equation}
where $\varepsilon_1>0$ will be chosen later (small enough).
As with $E$ and $M$, by a slight abuse of notation, we let $V_1(t)=V_1(y(t),y_t(t))=E(y(t),y_t(t)) + \varepsilon_1 M(y(t),y_t(t))$. Using \eqref{dotE-new} and \eqref{dotMleqnb3}, we have:
\begin{align}
\dot V_1 &= \dot E + \varepsilon_1 \dot M \nonumber \\
&\leq - \int_0^L a(x)y_t^2 \, dx - \varepsilon_1 \frac{\delta_p}{2} \int_0^L y_x^2 \, dx + \varepsilon_1 \int_0^L \left( p(x) - \frac{1}{2}q'(x) \right) y_t^2 \, dx \nonumber \\
&\quad + \varepsilon_1 C_1 \int_{x_0}^{x_1}  y^2 \, dx + \varepsilon_1 C_2(\delta_p/2) \int_0^L a(x) y_t^2 \, dx \text{ in } \mathcal{D}'(0,+\infty).
\label{dotV1leq}
\end{align}
Let us examine the term $\int \left( p - \frac{1}{2}q' \right) y_t^2 \, dx$ in \eqref{dotV1leq}.
From \eqref{prop:deltaq} and \eqref{support-p},
\begin{equation}
\label{pqoutside}
p - \frac{1}{2}q' \leq -\frac{\delta_q}{2}< 0 \text{ in } [0,x_0]\cup [x_1,L].
\end{equation}
Let
\begin{equation}
K := \max_{x \in [0,L]} \left( p(x) - \frac{1}{2}q'(x) \right)_+.
\end{equation}
We have, using \eqref{aona0a1}:
\begin{equation}
\label{majo-pqprime}
\int_{x_0}^{x_1} \left( p(x) - \frac{1}{2}q'(x) \right) y_t^2 \, dx \leq \frac{K}{a_{\text{min}}} \int_{x_0}^{x_1} a(x) y_t^2 \, dx \leq \frac{K}{a_{\text{min}}} \int_0^L a(x) y_t^2 \, dx.
\end{equation}
Let
\begin{equation}
\label{def:C3}
C_3 := \frac{K}{a_{\text{min}}} + C_2(\delta_p/2).
\end{equation}
From \eqref{dotV1leq}, \eqref{pqoutside}, \eqref{majo-pqprime}, and \eqref{def:C3},
\begin{multline}
\label{dotV1leqnb2}
\dot V_1 \leq - \left( 1 - \varepsilon_1 C_3 \right) \int_0^L a(x) y_t^2 \, dx
\\
- \varepsilon_1 \frac{\delta_p}{2} \int_0^L y_x^2 \, dx - \varepsilon_1 \frac{\delta_q}{2} \int_{[0,L]\setminus [x_0,x_1]} y_t^2 \, dx + \varepsilon_1 C_1 \int_{x_0}^{x_1} y^2 \, dx \text{ in } \mathcal{D}'(0,+\infty).
\end{multline}
We require
\begin{equation}\varepsilon_1 < \frac{1}{2C_3},
\end{equation}
so that, using also \eqref{a-continuous-and-geq0} and \eqref{aona0a1},
\begin{equation}
\label{1C3}
- \left( 1 - \varepsilon_1 C_3 \right) \int_0^L a(x) y_t^2 \, dx \leq - \frac{1}{2} \int_0^L a(x) y_t^2 \, dx \leq - \frac{a_{\text{min}}}{2} \int_{x_0}^{x_1} y_t^2 \, dx.
\end{equation}
Let us define
\begin{gather}
\label{def:delta2}
\delta_2 := \min\left\{\frac{a_{\text{min}}}{2}, \frac{\varepsilon_1 \delta_q}{2} \right\}> 0,
\\
\label{def:delta1}
\delta_1 := \frac{\varepsilon_1 \delta_p}{2}.
\end{gather}
From \eqref{dotV1leqnb2}, \eqref{1C3}, \eqref{def:delta2}, and \eqref{def:delta1}, we have
\begin{equation}
\label{eq:V1_bound}
\dot V_1 \leq - \delta_1 \int_0^L y_x^2 \, dx - \delta_2 \int_0^L y_t^2 \, dx + \varepsilon_1 C_1 \int_{x_0}^{x_1} y^2 \, dx \text{ in } \mathcal{D}'(0,+\infty).
\end{equation}
\begin{remark}
\label{rem:cross-terms}
\begin{enumerate}
\item The use of integral terms of the form $\int_0^L p(x) y y_t \, dx$ to study the behavior at infinity of solutions to the wave equation is now classical. To the best of our knowledge, the first article where this term is used to construct a strict Lyapunov function in the case $(x_0,x_1)=(0,L)$ is \cite{1979-Chen-SICON} (and even in the multi-dimensional case). For nonlinear wave equations, the first articles are, to the best of our knowledge, \cite{1978-Dafermos-Proceedings} (see, in particular, (3.15) in this article), \cite{1987-Haraux-book} (see, in particular, the definition of $\psi(t)$ on page~271), and \cite{1988-Haraux-Zuazua-ARMA,1988-Zuazua-AA}. For an example involving the linear wave equation with boundary conditions more complex than Dirichlet boundary conditions, see, for instance, \cite{2026-Dahmani-Chitour-Nguyen-Roman-JDE}.

\item The use of integral terms of the form $\int_0^L q(x) y_t y_x \, dx$ is also very classical. It comes from  the Rellich--Pohozaev  multiplier \cite{1940-Rellich-MZ, 1965-Pohozaev-SMD} (for elliptic equations) or the Morawetz multiplier \cite{1961-Morawetz-CPAM} (for wave equations). It was first used to prove observability results in \cite{1986-Ho-CRAS}, see also \cite{1988-Lions-book-vol-1}, and, in particular, Chapter 7 of this book.
\end{enumerate}
\end{remark}

\subsection{The nonlocal functional $E_0$}
\label{sec:nonlocal}
Looking at \eqref{eq:V1_bound}, we see that it remains to absorb the positive zeroth-order term $\int_{x_0}^{x_1} y^2 \, dx$. To do so, we are partly inspired by \cite{2025-Bastin-Coron-Hayat-Hal}, in particular Appendix~D of that article.
Let $\mathcal{L}: (f,g)\in H^1_0(0,L)\times L^2(0,L) \rightarrow h=\mathcal{L}(f,g)\in H^1_0(0,L)$ be the linear map defined by requiring
\begin{gather}
\label{def:L}
h_{xx}= a(x)f+g \text{ in }H^{-1}(0,L).
\end{gather}

Let us define, for $(f,g)\in H^1_0(0,L)\times L^2(0,L)$,
\begin{equation}
\label{def:E0}
E_0(f,g):= E(\mathcal{L}(f,g),f)\geq 0.
\end{equation}
Note that there exists $C_4>0$ such that
\begin{equation}
\label{majorantE0}
E_0(f,g)= E(\mathcal{L}(f,g),f)\leq C_4 \left(\|f\|_{L^2}^2+\|g\|_{H^{-1}}^2\right).
\end{equation}
With our usual abuse of notation, let
\begin{equation}
\label{def:E0(t)}
E_0(t):=E_0(y(t),y_t(t)), \quad t\geq 0.
\end{equation}
Let us check that $E_0\in C^1([0,+\infty))$ and that
\begin{equation}
\label{dotE0(t)}
\dot E_0(t)=-\int_0^L a(x)y(t,x)^2\, dx.
\end{equation}
Let $h:[0,+\infty)\times [0,L]\rightarrow \mathbb{R}$ be defined by
\begin{equation}
\label{def:h}
h(t,x):= \mathcal{L}(y(t),y_t(t))(x).
\end{equation}
As in the proof of \eqref{dotMleqnb1-dist}, we may assume that $y_0 \in H^2(0,L)\cap H^1_0(0,L)$ and $y_1 \in H^1_0(0,L)$. Then, according to Theorem~\ref{th:Cauchy-pb}, $y$ has the regularity given by~\eqref{eq-sol-reg}, so that
\begin{equation}
\label{regularity-h}
h\in C^0([0,+\infty);H^1_0(0,L)\cap W^{2,\infty}(0,L))\cap C^1([0,+\infty);H^1_0(0,L) \cap H^2(0,L)),
\end{equation}
and, using the definition of $\mathcal{L}$ and \eqref{def:h},
\begin{equation}
\label{eq:htxx}
(h_t)_{xx}=(h_{xx})_{t}=a(x)y_t+y_{tt}=y_{xx} \text{ in } C^0([0,+\infty);L^2(0,L)).
\end{equation}
Since both $h_t$ and $y$ vanish at $x=0$ and $x=L$, we deduce from \eqref{eq:htxx} that
\begin{equation}
\label{ht=y}
h_t=y.
\end{equation}
In particular,
\begin{equation}
\label{regularity-h-new}
h \in C^1([0,+\infty); H^1_0(0,L)\cap H^2(0,L)) \cap C^2([0,+\infty); H^1_0(0,L)) \cap C^3([0,+\infty); L^2(0,L)).
\end{equation}
From \eqref{def:h} and \eqref{ht=y}, one has
\begin{equation}
\label{edp:h}
h_{tt}=y_t=h_{xx}-a(x)y=h_{xx}-a(x)h_t.
\end{equation}
Since $h$ vanishes at $x=0$ and $x=L$, $h$ is a solution to our wave equation \eqref{eq-wave}, and we have
(see \eqref{dotE-new})
\begin{equation}
\label{deriv-enery-h}
\frac{d}{dt}E(h(t),h_t(t))=-\int_0^L a(x)h_t(t,x)^2 \, dx,
\end{equation}
which, together with \eqref{def:E0}, \eqref{def:E0(t)}, \eqref{def:h}, and \eqref{ht=y}, yields exactly \eqref{dotE0(t)}.

\subsection{Our strict Lyapunov function}
Finally, it suffices to define
\begin{equation}
\label{def:V}
V := V_1+C_0E_0,
\end{equation}
with
\begin{equation}
\label{defC0}
C_0:=\frac{\varepsilon_1 C_1}{a_{\text{min}}}.
\end{equation}
Note that, by \eqref{def:M}, \eqref{defV1}, \eqref{majorantE0}, and \eqref{def:V}, there exists $C_5>0$ such that
\begin{gather}
\label{Vlesenergy}
 V\leq C_5 E.
\end{gather}
Using the Young inequality, the optimal Poincar\'{e} inequality, and \eqref{def:M}, one has
\begin{equation}
\label{MleqE}
|M(f,g)|\leq \left(\|q\|_{L^\infty}+\frac{L}{\pi}\|p\|_{L^\infty}\right)E(f,g),\quad \forall (f,g)\in H^1_0(0,L)\times L^2(0,L).
\end{equation}
From now on let us assume that $\varepsilon_1$ is chosen so that
\begin{equation}
\label{varepsilon1<}
0<\varepsilon_1< \left(\|q\|_{L^\infty}+\frac{L}{\pi}\|p\|_{L^\infty}\right)^{-1}.
\end{equation}
Let
\begin{equation}
\label{def:deltaM}
\delta_M:= 1- \varepsilon_1 \left(\|q\|_{L^\infty}+\frac{L}{\pi}\|p\|_{L^\infty}\right)>0.
\end{equation}
Then, using in particular  \eqref{defV1}, \eqref{def:E0}, \eqref{def:V}, and \eqref{MleqE},
\begin{gather}
\label{VgeqE}
\delta_M E\leq V.
\end{gather}
Moreover, using \eqref{eq:V1_bound}, \eqref{dotE0(t)}, \eqref{def:V}, and \eqref{defC0}, one obtains
\begin{gather}
\label{dotVleq-gammaV}
\dot V \leq - \delta_1 \int_0^L y_x^2 \, dx - \delta_2 \int_0^L y_t^2 \, dx\leq -\gamma V \text{ in } \mathcal{D}'(0,+\infty),
\end{gather}
with
\begin{equation}
\label{def:gamma}
\gamma := 2\frac{\min\{\delta_1,\delta_2\}}{C_5}.
\end{equation}
Using also \eqref{VgeqE}, this concludes the proof that $V$ is a strict Lyapunov function. Finally, \eqref{Vlesenergy}, \eqref{VgeqE}, and \eqref{dotVleq-gammaV} also provide a proof of the exponential stability in $H^1_0(0,L)\times L^2(0,L)$, i.e., \eqref{exp-weak}, with $\gamma$ given by \eqref{def:gamma}.
\begin{remark}
Simple computations show that \eqref{Vlesenergy} holds for
\begin{equation}
\label{eq:possibleC5}
C_5=1+\varepsilon_1\left(\|q\|_{L^\infty}+\frac{L}{\pi}\|p\|_{L^\infty}\right)
+\frac{\varepsilon_1C_1}{a_{\text{min}}}\left(2\frac{L^4}{\pi^4}\|a\|_{L^\infty}^2+3\frac{L^2}{\pi^2}\right).
\end{equation}
It gives an explicit value for $\gamma$ defined in \eqref{def:gamma}.
\end{remark}

\section{Our localized strict Lyapunov function for the strong solutions}
\label{sec:strict-for-strong}
The Lyapunov function $V$ contains a ``nonlocal'' term (namely, $E_0$). However, if one considers solutions in $C^0([0,+\infty); \left(H^1_0(0,L)\cap H^2(0,L)\right)\times H^1_0(0,L))$ as in Theorem~\ref{th:Cauchy-pb}, one can take the Lyapunov function to be
\begin{equation}
V_{\text{s}}(y,y_t)=V_1(y_t,y_{tt})+C_0E(y,y_t)=V_1(y_t,y_{xx}-a(x)y_t)+C_0E(y,y_t),
\end{equation}
which now features only local terms. Note that $y_t$ is a weak solution to \eqref{eq-wave} (for the initial condition $(y_1, y_{0xx}-a(x) y_1)\in H^1_0(0,L)\times L^2(0,L)$) and that $h$ (see \eqref{def:h}) for this solution is simply $y$. Hence, applying \eqref{dotVleq-gammaV} to this weak solution, we obtain
\begin{gather}
\label{dotV1leq-new}
\dot V_{\text{s}} \leq  -\gamma V_{\text{s}} \text{ in } \mathcal{D}'(0,+\infty).
\end{gather}
One easily checks that there exists $C_6>0$ such that, for every $(y,y_t)\in \left(H^1_0(0,L)\cap H^2(0,L)\right)\times H^1_0(0,L)$,
\begin{gather}
\label{V1comparableH2H1}
\frac{1}{C_6} \left( \|y\|^2_{H^2}+\|y_t\|^2_{H^1}\right)  \leq V_{\text{s}} \leq C_6 \left( \|y\|^2_{H^2}+\|y_t\|^2_{H^1}\right).
\end{gather}
From \eqref{dotV1leq-new} and \eqref{V1comparableH2H1}, one obtains the exponential stability in $\left(H^1_0(0,L)\cap H^2(0,L)\right)\times H^1_0(0,L)$ (i.e., \eqref{exp-strong}) by a Lyapunov approach, with an exponential decay rate at least equal to $\gamma$ given by \eqref{def:gamma}.

\begin{remark}
Such a nonlocal/local phenomenon already appears in \cite[Section 6 and Appendix D]{2025-Bastin-Coron-Hayat-Hal}; see, in particular, Proposition~6 in that article, where a strict Lyapunov function for $L^2$ exponential stability is given by (D.2) and features two nonlocal terms (the first two terms of (D.2)), while the strict Lyapunov function for $H^1$ exponential stability given by (5.11) contains only local terms.
\end{remark}

\section{Conclusion and perspectives}
\label{sec:conclusion-and-perspectives}
In this paper, we have constructed an explicit, strict Lyapunov function for the one-dimensional linear wave equation subject to locally distributed damping and Dirichlet boundary conditions. A notable feature of our construction is that the Lyapunov function incorporates a nonlocal term for weak solutions, yet reduces to purely local terms in the case of strong solutions.

Our approach opens several natural directions for future research. For instance:
\begin{enumerate}
\item
Complex boundary conditions: Extending this approach to construct strict Lyapunov functions for more
intricate boundary configurations, such as those recently investigated in \cite{2024-Chitour-Nguyen-Roman-A}.
\item
Broader functional spaces: Generalizing the construction to more complex state spaces, particularly the $W^{1,p}_0(0,L) \times L^p(0,L)$
framework considered in \cite{2024-Chitour-Nguyen-COCV}.
\item Nonlinear one-dimensional wave equation: Extending this approach to construct strict Lyapunov functions in the nonlinear case. See, in particular,
\cite{2023-Vanspranghe-Ferrante-Prieur}.
\item
Multi-dimensional settings: Addressing the multi-dimensional wave equation with locally distributed damping.
A natural starting point would be a manifold without boundary, such as the two-dimensional torus $\mathbb{T}^2$:
\begin{multline}
y_{tt} - y_{x_1x_1} - y_{x_2x_2}
\\ + a(x)y_t = 0, \quad t \geq 0, \quad x=(x_1,x_2) \in \mathbb{T}^2.
\end{multline}
Inequality \eqref{aona0a1} is now replaced by the condition
\begin{equation}
a(x) \geq a_{\text{min}} > 0, \quad \text{for almost every } x \in \omega,
\end{equation}
where $\omega$ is a nonempty open subset of $\mathbb{T}^2$ satisfying the Geometric Control Condition
introduced by Bardos, Lebeau, and Rauch \cite{1992-Bardos-Lebeau-Rauch-SICON}.
 Their foundational work establishes that this condition guarantees exponential stability within the invariant space of states $(y,y_t)$ such that
\begin{equation}
\int_{\mathbb{T}^2} (y_t + a(x)y) \, dx = 0.
\end{equation}
An open question is whether one can  recover this theorem by means of an explicit, strict Lyapunov function.
\end{enumerate}

\section*{Acknowledgements}
The author thanks Yacine Chitour, Hoai-Minh Nguyen, and Christophe Prieur for useful comments, discussions, and references.

\printbibliography

\end{document}